
\documentstyle{amsppt}
\pageheight{46pc}
\pagewidth{33pc}
\redefine\P{\bold P}   
\define\Q{\bold Q}     
\define\C{\bold C}     
\define\Z{\bold Z}     
\define\R{\bold R}          
%

%


\define\cI{\Cal I}      


\define\cX{\Cal X}
\define\cY{\Cal Y}
\define\Sym{\operatorname{Sym}}     

\define\codim{\operatorname{codim}}
\define\gr{\operatorname{gr}}
\define\Spec{\operatorname{Spec}}
\define\CH{\operatorname{CH}}
\redefine\H{\operatorname{H}}
\define\N{\operatorname{N}}

\define\ie{i\.e\. }
\define\eg{e\.g\. }

\define\isom{\cong}
\define\Oh#1{{\Cal O}_{#1}}
\define\ext#1{\overset{#1}\to{\Lambda}}
\redefine\mod#1{\mid\!#1\!\mid}
\define\tensor{\otimes}
\define\onto{\hbox{$\to\!\!\!\to$}}
\define\into{\hookrightarrow}
\define\from{\leftarrow}
\define\contain{\supset}
\define\propersubset{\mathrel{\lower4pt\hbox{$\scriptscriptstyle\not$}%
\mkern-5mu\subseteq}}
\define\restrict#1{\!\mid_{#1}}
\magnification\magstep1
\NoRunningHeads
\topmatter
\title Cohomological and Cycle-theoretic Connectivity \endtitle
\author Kapil H. Paranjape \endauthor
\affil Tata Institute of Fundamental Research, Bombay. \endaffil
\address School of Mathematics,
        Tata Institute of Fundamental Research,
        Homi Bhabha Road,
        Bombay 400005,
        INDIA \endaddress
\email KAPIL\@TIFRVAX.bitnet \endemail
\date 4 April, 1992 \enddate
\thanks The author would like to thank the University of Chicago for
support during the period October 1990--September 1991. \endthanks
\abstract One of the themes in algebraic geometry is the study of the
relation between the ``topology'' of a smooth projective variety and a
(``general'') hyperplane section.  Recent results of Nori produce
cohomological evidence for a conjecture that a general hypersurface of
sufficently large degree should have no ``interesting'' cycles. We
compute precise bounds for these results and show by example that
there are indeed interesting cycles for degrees that are not high
enough. In a different direction Esnault, Nori and Srinivas have shown
connectivity for intersections of small multidegree. We show analogous
cycle-theoretic connectivity results.
\endabstract
\endtopmatter
\document

\heading Introduction \endheading

We study the relation between the ``topology'' of a smooth projective
variety and a ``general'' subvariety. One of the measures of topology
is a suitable cohomology theory, another measure is the group of
cycles modulo rational equivalence---the Chow group. These two are
conjecturally related by series of conjectures of A. A. Beilinson and
S. Bloch (see \cite{11}). We formulate and examine some concrete
cases.

The classical ``weak'' Lefschetz theorem states that if $X$ is a
smooth projective variety and $Y$ an ample divisor then
$\H^l(X)\to\H^l(Y)$ is an isomorphism for $l<\dim Y$. The conjectural
cycle-theoretic analogue is that $\CH^p(X)_{\Q}\to\CH^p(Y)_{\Q}$ is an
isomorphism for $p<\dim Y/2$ (see (1.5)). The case $p=1$ is a
classical theorem of S. Lefschetz and A. Grothendieck (see \cite{9}).
We prove the conjecture in some cases (a more precise statement is
(5.4)).
\proclaim{Theorem 0.1} Given integers $1\leq d_1\leq\ldots\leq d_r$
and any non-negative integer $l$, let $X\subset\P^n$ be a smooth
subvariety of multidegree $(d_1,\ldots,d_r)$. If $n$ is sufficiently
large then $\CH^l(X)_{\Q}\cong\Q$. \endproclaim

The classical theorem of M. Noether and S. Lefschetz states that if
$X=\P^3$ and $Y\subset X$ is a ``general'' surface with $\deg Y\geq
4$, then $\CH^1(Y)=\Z$.  There are also some recent generalisations by
M. Green \cite{6}. In a recent paper \cite{14} M. V. Nori shows that
we have a cohomological analogue of this statement as well.  In
particular, from the conjectural framework mentioned above one should
expect that if $Y$ is a ``general'' complete intersection in $X$ of
{\it sufficiently high multidegree} then $\CH^p(X)\to\CH^p(Y)$ is an
isomorphism for $p<\dim Y$ (see (2.9)). We compute precise bounds for
the degrees required. As a result we find that if $X=\P^n$ then we
need $\deg Y\geq (2n -2)$ in Nori's result (see Section~2). To
contrast this we show (see also \cite{18}):
\proclaim{Proposition 0.2} Let $Y$ be a ``general'' hypersurface in
$\P^n$ for $n\geq 3$ of degree between $n$ and $2n-3$. Then $Y$
contains a pair of lines $L_1$ and $L_2$ such that the difference is
not rationally equivalent to zero in the Chow group of 1-cycles on
$Y$. \endproclaim

It was shown by A.~A.~Roitman \cite{16} that the Chow group of
0-cycles on intersections of small multidegree is $\Z$.  Recent
results of H. Esnault, M. V. Nori and V. Srinivas \cite{4} give a
cohomological analogue for schemes $Y\subset\P^n$ defined by a small
number of equations of low degree. These results led Srinivas to
suggest that there are other ``rational-like'' connectivity properties
of such varieties; \ie we should have $\CH_p(Y)_{\Q}=\Q$ for small
enough $p$ (see (1.9)). The precise prediction for cubic hypersurfaces
$Y$ is that if $\dim Y>4$ then $\CH_1(Y)_{\Q}=\Q$. We prove this for
the general cubic hypersurface in (4.1.4). More generally, we prove
the following theorem (a more precise formulation of the result is
(4.2.5)).
\proclaim{Theorem 0.3} Given integers $1\leq d_1\leq\ldots\leq d_r$
and any non-negative integer $l$, choose any $r$ sections
$f_i\in\Gamma(\P^n,\Oh{\P^n}(d_i))$ and let $X=V(f_1,\ldots,f_r)$ be
the scheme defined by these. If $n$ is sufficiently large then
$\CH_l(X)\cong\Z$ and is generated by the class of a linear subspace
$\P^l\subset X$. \endproclaim

In Section~1 we define the notion of cohomological and cycle-theoretic
connectivity. We also state the general conjectures that form a basis
for the kind of results one is looking for. We follow up in Section~2
with a sketch of the proof of the Theorem of Nori where we also obtain
the bounds for the numbers $N(e)$ appearing in this Theorem. We
produce examples in Section~3 to show the sharpness of these bounds.
In Section~4 we study cycles of small dimension on intersections of
small multidegree in projective space. Finally, in Section~5 we obtain
results on cycles of small codimension on such varieties.

I am very grateful to M. V. Nori for explaining his ideas and results
related to the papers \cite{14} and \cite{4}. C. Schoen gave me his
preprint \cite{18} on connectivity and also clarified many ideas on
cycles. The question raised therein on cubic hypersurfaces led to
Section~(4.1) of this paper. Discussions with K. Joshi, B. Kahn, N.
Mohan Kumar and V. Srinivas were very helpful. In particular, B. Kahn
pointed out the paper of Leep and Schmidt \cite{12}. V. Srinivas
helped me solve the optimisation problem that appears in Section~2 and
provided much guidance and encouragement.

\heading 1. Preliminaries on Connectivity \endheading

We work over a fixed base field of characteristic zero (say $\C$)
which we denote by $k$ and a fixed universal domain $K\contain k$; \ie
$K$ is an algebraically closed field of infinite transcendence degree
over $k$. Let $X$ be a scheme of finite type over $K$. Let $\H^*(X)$
denote the de~Rham cohomology of $X$ {\it relative to $k$} (see
\cite{10}). This is the {\it direct limit} of the de~Rham cohomology
groups $\H^*(\cX)$ as $\cX$ runs over all models of $X$ which are of
finite type over $k$.

\definition{Definition 1.1} Let $X$ be a variety over $K$ and
$Y\subset X$ a closed subvariety. We say that the inclusion of $Y$ in
$X$ is a {\it cohomological $r$-equivalence} if the restriction
morphism $\H^l(X)\to\H^l(Y)$ is an isomorphism for $l\leq r$.
\enddefinition

\definition{Definition 1.2} Let $X$ be a variety over $K$ and
$Y\subset X$ a closed subvariety. We say that the inclusion of $Y$ in
$X$ is a {\it cycle-theoretic $c$-equivalence} if the restriction
morphism $\CH^p(X)_{\Q}\to\CH^p(Y)_{\Q}$ is an isomorphism for
$p\leq c$. \enddefinition

A typical situation is the following. We have a smooth projective
variety $X$ {\it over $k$} and $\cY\subset X\times S$ is a family of
smooth subvarieties. We choose a $k$-embedding of $k(S)$ into $K$,
where $k(S)$ is the function field of $S$, and put $Y=\cY\times_S\Spec
K$.  This corresponds to the study of restriction of cycles (or
cohomology classes) to the ``general'' subvariety of $X$ in the family
parametrised by $S$.

\remark{Remark 1.3} The classical theorem of S. Lefschetz (the
``weak'' Lefschetz theorem) says that if $Y$ is an ample divisor in a
smooth projective variety $X$, then the inclusion $Y\into X$ is a
cohomological $r$-equivalence, for $r=\dim Y-1$. By induction we have
the same result for a complete intersection subvariety.\endremark

\remark{Remark 1.4} One expects the Chow group $\CH^p(X)_{\Q}$ to
carry a natural filtration $F^{\cdot}$ such that we have a relation
between $\gr^l_F\CH^p(X)_{\Q}$ and $\H^{2p-l}(X)$ (see \cite{11}).
More precisely, we may expect cohomological $r$-equivalence to imply
(tensor $\Q$) cycle-theoretic $c$-equivalence where $2c\leq r$. In
view of (1.3) there is some interest in determining cycle-theoretic
connectivity for ample divisors (and more generally complete
intersection subvarieties). \endremark

With notation as in (1.3), the theorems of S.~Lefschetz and
A. Grothendieck (see \eg \cite{9}) say that the inclusion $Y\into X$
is a cycle-theoretic 1-equivalence if $\dim Y\geq 3$.  This is a
particular case of the following:
\proclaim{Conjecture1.5} Let $Y$ be a smooth ample divisor in $X$,
then the inclusion of $Y$ in $X$ is a cycle-theoretic $c$-equivalence
for $2c\leq\dim Y-1$.\endproclaim
We obtain some results in this direction in Section~5.

\remark{Remark 1.6} The theorem of M. Noether and S. Lefschetz says
that if we take $Y$ to be the ``general'' complete intersection in
$X=\P^n$ then the inclusion of $Y$ in $X$ is a cycle-theoretic
1-equivalence if $\dim Y\geq 2$. The following theorem of Nori is a
cohomological generalisation:
\proclaim{Theorem}\rom{(M. V. Nori)} Let $X$ be a smooth projective
variety over $k$ and let $\Oh{X}(1)$ be an ample line bundle on
$X$. Let $V\to\Gamma(X,\Oh{X}(1))$ be a space of sections which
generates $\Oh{X}(1)$ at stalks. Let $d_1\leq\cdots\leq d_r$ be a
multidegree and $Y$ be the ``general'' complete intersection in $X_K$
of this multidegree. Then the inclusion $Y\into X_K$ is a
cohomological $(\dim Y+e-1)$-equivalence, if $d_1\geq N(e)$, where
$N(e)$ is an integer depending only on $(X,\Oh{X}(1),V)$ and
$e\leq\dim Y$. \endproclaim
By (1.4) this leads to the following conjecture also due to Nori:
\proclaim{Conjecture}\rom{(M. V. Nori)} In the situation of the above
theorem, the inclusion of $Y$ in $X_K$ is a cycle-theoretic
$c$-equivalence for $2c\leq\dim Y +e-1$.\endproclaim
Note that we are free to choose $e$; however, if we want an
isomorphism $\CH^l(X_K)\isom\CH^l(Y)$ for all $l\leq (\dim Y+e-1)/2$
as predicted, then the best one can hope for is $e=\dim Y -1$ by the
results of D.~Mumford and A.~A.~Roitman \cite{17} (see also
Example~(3.1)). \endremark

We will make this conjecture more precise in Section~2 by computing a
value for $N(e)$. In Section~3 we will show that for $X=\P^n$ this
gives a sharp bound for cycle-theoretic connectivity.

\remark{Remark 1.7} Reverting back to varieties over $k$ we know that
the de~Rham cohomology theory satisfies the Hard Lefschetz theorem for
smooth projective varieties $X$ over $k$. As shown in \cite{11} a
consequence of the conjectured relation mentioned in (1.3) and the
Hard Lefschetz theorem is that $\gr^l_F\CH^p(X)_{\Q}$ should depend
only on $\H^{2p-l}(X)/\N^{p-l+1}\H^{2p-l}(X)$ for $l\leq p$, where
$\N^{\cdot}$ denotes the filtration by coniveau (see \cite{8}).
\endremark

There is a natural filtration on de~Rham cohomology---the Hodge
filtration---which conjecturally yields the coniveau filtration (via
the comparison with singular cohomology; see \cite{8}). With this in
mind we examine a result of Esnault-Nori-Srinivas \cite{4}:
\proclaim{Theorem}\rom{(H. Esnault, M. V. Nori, V. Srinivas)}
Let $X=V(f_1,\ldots,f_r)\subset\P^n$ be the zero locus of a
collection of homogeneous equations $f_i$ of degrees $d_i$, where we
assume $d_1\leq\cdots\leq d_r$. Then for any $i\geq 0$,
$$ F^k\H^i_c(\P^n-X) = \H^i_c(\P^n-X),$$
where $F^{\cdot}$ denotes the Hodge filtration and
$$ k = \left[ {{n-\sum_{j=2}^r d_r}\over{d_1}} \right].$$
\endproclaim
This leads us to the following:
\proclaim{Conjecture 1.9} If $X$ is a smooth subvariety of
$\P^n$ of multidegree $d_1\leq\cdots\leq d_r$, then
$\CH_p(X)_{\Q}=\Q$ for $p<k$ with $k$ as above.\endproclaim

\remark{Remark 1.10} For $n\geq\sum_{j=1}^r d_r$, this conjecture
says that $\CH_0(X)_\Q=\Q$. This particular case has been proved by
A. A. Roitman \cite{16} even with $\Z$ instead of $\Q$. We will
obtain generalisations in Section~4. \endremark

\heading 2. Cohomological connectivity \endheading

In this section we study cohomological connectivity for the
``general'' member of the family of all complete intersection
subvarieties of a fixed multidegree in a smooth projective variety.
We sketch a proof of the theorem of Nori stated in Section~1. In
addition, we will obtain bounds (see (2.8)) for the number $N(e)$
apppearing in this theorem in terms of some natural invariants for the
triple $(X,\Oh{X}(1),V)$.

A more general situation than the one in the Theorem is the following.
Let $X$ be a smooth projective variety and $F$ a vector bundle on $X$
which is generated at stalks by a space of sections $W\to\Gamma(X,F)$.
We define the vector bundle $M$ by the exact sequence
$$ 0 \to M \to W \tensor \Oh{X} \to F \to 0 \tag 2.2$$
Let $S$ be the projective space $\P(W^*)$. We have $\cY=\P_X(M^*)$
which is a closed subvariety of $X\times S$. The morphism $\cY\to S$
makes it the family of subvarieties of $X$ defined by zeroes of
sections of $E$. Choose a $k$-embedding into $K$ of the function
field of $S$ and let $Y=\cY\times_S\Spec K$. Let the sheaves
$\Omega^i_{(X\times S,\cY)}$ be defined by the exact sequence
$$ 0 \to \Omega^i_{(X\times S,\cY)} \to \Omega^i_{X\times S}
                \to \Omega^i_{\cY} \to 0 \tag 2.3$$
Madhav Nori \cite{14} has shown that his Theorem follows from
the vanishing of the direct image sheaves
$\R^a {p_S}_*(\Omega^b_{(X\times S,\cY)})$ for $a\leq\dim Y$
and $a+b\leq\dim Y+e$.

\remark{Remark 2.4} V. Srinivas \cite{19} has defined Hodge cohomology
for any scheme over $k$ as follows:
$$\H^n(X)=\oplus_{p=0}^n\H^{n-p}(X,\Omega^p_{X/k}).$$
He has constructed Gysin and cycle class maps for this theory when $X$
and $Y$ are smooth varieties over {\it any} field extension $K$ of
$k$. If this is a Poincar\'e duality theory, then one may view the
result of M.~V.~Nori stated above as a result on Hodge-theoretic
connectivity.  \endremark

\proclaim{Lemma 2.5} The higher direct image sheaf
$\R^a{p_S}_*(\Omega^b_{(X\times S,\cY)})$ vanishes if the cohomology
groups $H^{k}(X,\ext{m}{M}\tensor\ext{l}{F}\tensor\Omega^c_X)$ vanish
for all non-negative integers $k$, $l$, $m$ and $c$ such that
$k\leq\dim X-a$, $m+(\dim X -c)\leq b$ and $k+l+m-c = b+1-a$.
\endproclaim

\demo{Sketch of Proof}
As in \cite{14} an $E_1$ spectral sequence argument
shows that $\R^a {p_S}_*(\Omega^b_{(X\times S,\cY)})$ vanishes if the
sheaves
$\R^a {p_S}_*(\Omega^{\beta}_{(X\times S,\cY)/X}
                \tensor\Omega^{\gamma}_X)$
vanish for $\beta+\gamma=b$; where the sheaf
$\Omega^{\beta}_{(X\times S,\cY)/X}$ is defined by an exact sequence
analogous to (2.3) with relative differentials. These sheaves can be
computed by using a natural resolution for the sheaf
$\Omega^b_{(X\times S,\cY)/X}$ which we now explain.

Let $G_r$ be the filtration of $\ext{k}{W^*}$ obtained by taking
exterior powers of the exact sequence (2.2); the numbering is so
chosen that $\gr_r^G(\ext{k}{W^*})=\ext{k-r}{F^*}\tensor\ext{r}{M^*}$.
Let $\Oh{S}(1)$ be the tautological line bundle. We use the notation
$G_b(\ext{k}{W^*})(-k)$ to denote the sheaf $p_X^*
G_b(\ext{k}{W^*})\tensor p_S^*\Oh{S}(-k)$ on $X\times S$.  We have
natural homomorphisms of sheaves on $S$, for each $k$ and $b$,
$$ G_b(\ext{k+1}{W^*})(-k-1) \to
                G_b(\ext{k}{W^*})(-k)    $$
which fit together for $k\geq b+1$ and give a resolution for
$\Omega^b_{(X\times S,\cY)/X}$, \ie we have an exact sequence.
$$ 0 \from \Omega^b_{\cY/X}
        \from \Omega^b_{X\times S/S}
         \from G_b(\ext{b+1}{W^*})(-b-1)
          \from G_b(\ext{b+2}{W^*})(-b-2) \cdots        $$
Thus the sheaf
$\R^a {p_S}_*(\Omega^b_{(X\times S,\cY)/X}\tensor\Omega^c_X)$
is the zero sheaf if for all
$\alpha\geq a$ and $\beta=b+1+(\alpha-a)$, the sheaves
$\R^{\alpha} {p_S}_*(G_b(\ext{\beta}{W^*})(-\beta)\tensor\Omega^c_X)$
are zero. By the projection formula we see that this sheaf is a twist
of the trivial sheaf on $S$ with fibre
$\H^{\alpha}(X,G_b(\ext{\beta}{W^*})\tensor\Omega^c_X)$
by $\Oh{S}(-\beta)$. Now, by the $E_1$ spectral sequence associated
with the filtration $G_{\cdot}$, we see that the group
$\H^{a}(X,G_l(\ext{b}{W^*})\tensor\Omega^c_X)$
vanishes if the groups
$\H^{a}(X,\ext{d}{M^*}\tensor\ext{b-d}{F^*}\tensor\Omega^c_X)$
vanish for all $d\leq l$. Applying Serre's duality theorem to
the latter groups we obtain
$\H^{n-a}(X,\ext{d}{M}\tensor\ext{b-d}{F}\tensor\Omega^{n-c}_X)$,
where $n=\dim X$. Combining the equations and inequalities above we
have the result.\qed\enddemo

Let us now revert to the earlier situation where we are given
$(X,\Oh{X}(1),V)$ as in the Theorem of Nori. Assume we are given
$d_1\leq\cdots\leq d_r$, a multidegree and let
$F=\oplus_i\Oh{X}(d_i)$.  Let us choose $W$ to be the direct sum of
$\Sym^{d_i}(V)$ over all $i$. If $M_k$ is defined by the exact
sequence
$$ 0 \to M_k \to \Sym^k(V)\tensor\Oh{X} \to \Oh{X}(k) \to 0, $$
then $M=\oplus_i M_{d_i}$ and so we can apply the following:

\proclaim{Lemma 2.6} Let $m_b$ be such that
$H^a(X,\Omega_X^b\tensor\Oh{X}(m_b-a))=0$ for all $a\geq 1$. Then
$$ H^a(X,\Omega_X^b\tensor\ext{c}{M}\tensor\Oh{X}(k)) = 0$$
for all $a\geq 1$ such that $a + k \geq m_b + c$.
\endproclaim

\demo{Proof} It is well known that the given condition on $\Omega_X^b$
is the condition for $m_b$-regularity in the sense of Castelnouvo and
Mumford (see \eg \cite{13}), which in turn implies that
$$ H^a(X, \Omega_X^b\tensor\Oh{X}(k)) = 0 $$
for all $a\geq 1$ such that $a + k \geq m_b$. Let $f: X \to \P(V)$
be the morphism given by $(V,\Oh{X}(1))$. We have a vector bundle
$E_k$ on $\P(V)$ defined by the exact sequence
$$ 0 \to E_k \to \Sym^k(V) \tensor \Oh{\P(V)}
                                \to \Oh{\P(V)}(k) \to 0 .$$
Clearly, we have $M_k = f^{*}E_k$. Moreover, $f$ is finite and so we
may compute cohomologies of sheaves on $X$ via their direct images
under $f$. In particular, we see that $f_{*}(\Omega_X^b)$ is
$m_b$-regular on $\P(V)$. By the projection formula we have
$$
f_{*}(\Omega_X^b\tensor\ext{c}{M}\tensor \Oh{X}(k)) =
        f_{*}\Omega_X^b\tensor\ext{c}{(\oplus_i E_{d_i})}
                        \tensor \Oh{\P(V)}(k) .$$
It is well known (see \eg \cite{5}) that $E_k$ is 1-regular. A direct
sum of 1-regular sheaves is also 1-regular.  Moreover, the regularity
of the tensor product of a coherent sheaf with a vector bundle on
$\P(V)$ is the sum of their regularities (see \cite{5}). Thus we have
the result.  \qed\enddemo

By this lemma we only need $k+d^{(l)} \geq m+m_c$ and $k\geq 1$ for
all integers $k$, $l$, $m$ and $c$ occurring in Lemma~(2.5); here
$d^{(l)}=\sum_{i=1}^l d_l$. An easy computation shows that $k\geq 1$
and $l\geq 1$ for all tuples $(k,l,m,c)$ satisfying the conditions in
Lemma~(2.5). Thus we need only take maximum of $(m+m_c-k)/l$ over all
such tuples as $(a,b)$ vary over all pairs such that $a\leq\dim Y$ and
$a+b\leq\dim Y+e$. This maximum will serve as $N(e)$.

For any triple $(X,\Oh{X}(1),V)$ we define
$$ m_X=\max\{m_c -c -1 : 0 \leq c \leq \dim X\}. \tag 2.7$$
One can solve the above optimisation problem after substituting $m_c$
by $m_X+c+1$ ; the maximum obtained this way is
$$ N(e) = \dim Y + e+1 + m_X. \tag 2.8$$
Thus we may restate the conjecture of M. V. Nori as follows:
\proclaim{Conjecture 2.9} Let $Y$ be the general complete intersection
of multidegree $d_1\leq\cdots\leq d_r$ in $X$ and assume that
$d_1\geq\dim Y+e+1+m_X$ for some $e\leq\dim Y-1$. Then the inclusion of
$Y$ in $X_K$ is a cycle-theoretic $c$-equivalence for $2c<\dim Y+e$.
\endproclaim

\remark{Remark 2.10} The number $m_c$ does not seem to have been
computed for general triples $(X,\Oh{X}(1),V)$, even when one assumes
reasonable properties. Work of L. Ein and R. Lazarsfeld \cite{3} shows
that $m_0 \leq (d_1 + \cdots d_c -c +1)$ where $c= \codim_{\P(V)} X$.

At the other extreme, by the Kodaira Vanishing  Theorem we have
$m_n \leq n+1$ for $n= \dim X$.
\endremark

\remark{Remark 2.11} For $X=\P^n$, we see easily that $m_c= c+1$ for
all $c$ in the range $0\leq c \leq n$; thus $m_{\P^n}=0$. Thus we have
a special case of (2.9):
\proclaim{Conjecture 2.12} Let $Y$ be the general hypersurface in
$\P^n$ of degree $\geq 2n-2$, then $\CH^p(Y)_{\Q}=\Q$ for $p<\dim Y$.
\endproclaim
In Section~3 we will show that this bound is sharp by showing that the
general hypersurface of degree $2n -3$ does in fact contain
interesting cycles. \endremark

\heading 3. Examples of large degree \endheading

In this section we work over a universal domain (say $\C$). We will
construct examples of hypersurfaces $X$ of degree $\leq 2n-3$ in
$\P^n$ such that $\CH_k(X)$ is not $\Z$ for $k=0,1$; \ie
cycle-theoretic connectivity in codimension $(\dim X -k)$ does not
hold.

\proclaim{Example 3.1} Let $X$ be a smooth projective variety and
$\{p_1,\ldots,p_n\} =S$ be a set of distinct points on $X$. There is a
divisor $Y\subset X$ such that $S\subset Y$ and
$$  \oplus_{i=1}^n\Z\cdot p_i \into \CH_0(Y), $$
\ie the points are linearly
independent in the group of 0-cycles modulo rational equivalence on
$Y$. In fact, if $\Oh{X}(1)$ is ample and $d$ sufficiently
large, then a ``general'' element $Y$ of the complete linear system
$\mod{\Oh{X}(d)}$ will have this property.  \endproclaim

\demo{Proof} Let $A$ be any line bundle on $X$ with the property that
there is a base point free linear system $V \subset \Gamma(X,Z)$ such
that the evaluation map $ V\tensor\Oh{X} \onto \oplus_{i=1}^n A_{p_i}
$ surjects. Then we have a morphism $f: X \to \P(V)$ such that $f(S)$
consists of linearly independent points in $\P(V)$.

\proclaim{Sublemma 3.2} If $H\subset\P(V)$ is a general hypersurface
of degree $\geq \dim V$ which contains $f(S)$ then $$ \oplus_{i=1}^n
\Z\cdot f(p_i) \into \CH_0(H) .$$ \endproclaim

Assuming this, we see that is $s\in\Sym^k(V)$ is a general element and
$k\geq \dim V$, then the zero locus $Y$ of $s$ in $X$ is a smooth
divisor. Now from the commutative diagram
$$
    \matrix\format\r&\quad\c&\quad\l\\
    \oplus_{i=1}^n \Z\cdot{p_i} & {\rightarrow} & \CH_0(Y) \\
    {\downarrow f_{*}} & & {\downarrow f_{*}}\\
    \oplus_{i=1}^n \Z\cdot{f(p_i)} & {\hookrightarrow} & \CH_0(H)
    \endmatrix
$$
we see that we have the result.
\qed\enddemo

\demo{Proof}(of (3.2)) Let $H'$ be a general hypersurface in $\P(V)$
of degree $\geq \dim V$. Then by the theorem of Mumford and Roitman
(see \eg \cite{17}), any $n$ general points $q_1,\ldots,q_n$ on $H'$
are linearly independent in the Chow group of 0-cycles on $H'$.

Now the set $f(S)$ consists of $n$ linearly independent points in
$\P(V)$. Hence there is an automorphism of $\P(V)$ that takes this
set into the set $q_1,\ldots,q_n$ (which by generality may be
assumed to be linearly independent also). Under this automorphism we
pull back the hypersurface $H'$ to give us the required hypersurface
$H$.
\qed\enddemo

\proclaim{Example 3.3} Let $L_1$, $L_2$ be a pair of skew lines in
$\P^n$ and $X$ a general hypersurface of degree $\geq n$ containing
$L_1$ and $L_2$.  Then we have
$$  \Z\cdot L_1 \oplus \Z\cdot L_2 \into CH_1(X) .$$
\endproclaim

\demo{Proof} Let $\P^{n-1}$ be a general linear space and
$p_i = L_i\cap \P^{n-1}$. We have seen during the proof of the
earlier example that if $Y$ is a general hypersurface of degree
$\geq n$ in $\P^{n-1}$ containing $p_1$ and $p_2$, then
$$ \Z\cdot p_1 \oplus \Z\cdot p_2 \into CH_0(Y) .$$

We have a short exact sequence
$$ 0 \to \cI_{L_1\cup L_2 /\P^n}(d-1) \to \cI_{L_1\cup L_2 /\P^n}(d)
        \to \cI_{p_1\cup p_2 /\P^{n-1}}(d) \to 0 .$$
Thus, the general hypersurface of degree $d$ in $\P^n$ containing
$L_i$ restricts to the general hypersurface of the same degree in
$\P^{n-1}$ containing $p_i$, providing
$H^1(\P^n,\cI_{L_1\cup L_2 /\P^n}(d-1))$ vanishes. Now we use the
exact sequence
$$ 0 \to \cI_{L_1\cup L_2 /\P^n}(d) \to \Oh{\P^n}(d) \to
        \Oh{L_1}(d) \oplus \Oh{L_2}(d) \to 0 $$
to see that $H^1(\P^n,\cI_{L_1\cup L_2 /\P^n}(d))=0$ for all $d\geq
1$. We have a diagram
$$
    \matrix\format\r&\c&\l&\quad\c&\quad\l\\
    \Z\cdot{L_1}&\oplus&\Z\cdot{L_2} & {\rightarrow} & \CH_0(X) \\
    &{\downarrow i^{*}}& & & {\downarrow i^{*}}\\
    \Z\cdot{p_1}&\oplus&\Z\cdot{p_2} & {\hookrightarrow} & \CH_0(Y)
    \endmatrix
$$
where $i: \P^{n-1} \to \P^n$ is the inclusion. Hence the result.
\qed\enddemo

\remark{Remark 3.4} By general considerations (see \eg \cite{15}) we
can show that a hypersurface in $\P^n$ of degree $\leq 2n-3$ always
contains a pair of lines. Hence we have Proposition~(0.2). \endremark

\heading 4. Examples of small degree \endheading

Fix a multidegree $d_1\leq\cdots\leq d_r$. If we take a large $n$ and
look at the varieties $X$ given by a collections of homogeneous
equations of the given multidegree in $\P^n$, then we get connectivity
statements for $\CH_k(X)$ for $k$ small. In the first part we study
the Chow group of 1-cycles on cubics as an illustrative example. In
the second part we prove the more general assertion (0.2).

\subheading{4.1 Cubic hypersurfaces}

We illustrate the general case by deducing that the Chow group of
1-cycles on a general cubic hypersurface of dimension 5 (over the
field of complex numbers) is $\Z$, using as starting point the
following well known facts.

\remark{Fact 4.1.1} The Chow group of 1-cycles on a quadric
hypersurface $X\subset \P^n_k$ of dimension at least 2 (\ie $n \geq
3$) is generated by lines, providing $k$ is algebraically closed.
\endremark

\remark{Fact 4.1.2} A quadric hypersurface $X\subset \P^n_k$ of
dimension at least 1 (\ie $n \geq 2$) contains a $k$-rational point
$p$ provided that $k$ is a $C_1$-field (in the sense of Lang; see \eg
\cite{7}), and hence $\CH_0(X)\cong \Z\cdot p$. \endremark

By general principles (see \eg \cite{15}) we see that if $k$ is an
algebraically closed field, a cubic hypersurface $X\subset\P_k^6$
contains a plane $P$. We project from this plane to get a morphism
$\tilde{X} \to \P_k^3$, where $\tilde{X}\to X$ is the blow up of $X$
along the plane. Let $E\subset\tilde{X}$ be the exceptional divisor of
this blow up. Then we see that $E\to\P_k^3$ is a family of conics and
$\tilde{X}\to\P_k^3$ is a family of quadric surfaces.

Let $C$ be any curve in $\tilde{X}$. If its image $D$ in $\P_k^3$ is a
curve, then the field $k(D)$ is $C_1$. The conic $E_{k(D)}=E
\times_{\P_k^3} \Spec k(D)$ then has a $k(D)$-rational point and so we
have a morphism $\tilde{D}\to E$ over $D$, where $\tilde{D}$ is the
normalisation of $D$. The 1-cycle $\xi = C - \deg(C/D)\cdot \tilde{D}$
restricts to a 0-cycle of degree 0 on the fibre of $\tilde{X}$ over
the point $\Spec{k(D)}\to \P_k^3$. By Fact~1.2, we see that this class is
then rationally trivial in this fibre. We can then apply the following
well known lemma to the cycle $\xi$ on $\tilde{X}\times_{\P_k^3} D$.

\proclaim{Lemma 4.1.3} Let $X\to Y$ be a proper morphism where $Y$ is
integral. Let $K$ be the function field of $Y$ and $X_K=X\times_Y
\Spec K$. If $\xi\in\CH_k(X)$ is such that it restricts to zero in
$\CH_{k'}(X_K)$, where $k'=k-\dim Y$, then there is a proper subscheme
$Z\propersubset Y$ such that $\xi$ is supported on $X\times_Y Z$.
\endproclaim

Thus $\xi$ is rationally equivalent on $\tilde{X}$ to a cycle
supported on the fibres of $\tilde{X}\to\P_k^3$.  By Fact~1.1 such
1-cyles are generated by lines. Hence we have a surjection
$$ \CH_1(E) \oplus \{\hbox{subgroup  generated by lines}\} \onto
        \CH_1(\tilde{X}).$$
But the morphism $\CH_1(E)\to\CH_1(X)$ factors through $\CH_1(P)$ which
is also generated by lines. Thus we see that the Chow group of 1-cycles
on $X$ is generated by lines. An easy computation shows that the
variety of lines on such a cubic is a Fano variety. If this variety is
smooth then it is known (see \cite{2}) that all lines on $X$ are
rationally equivalent. A similar argument will also give the following
result for $X$ of dimension bigger than 5. An analogous statement
appears in the paper of Schoen (see Theorem~5.1 in \cite{18}) for
cubics of dimension at least 6.

\proclaim{Proposition 4.1.4} Let $X\subset\P_k^n$ be a smooth cubic
hypersurface of dimension at least 5, such that the variety of lines
on it is smooth; then $\CH_1(X)\cong\Z$. \endproclaim

\subheading{4.2 Intersections of small multidegree}

Let $k$ be any field. Given positive integers $d_1\leq\cdots\leq d_r$,
and a non-negative integer $m$, we define $n=n(d_1,\ldots,d_r;m;k)$ to
be the smallest integer (possibly infinite!) such that for all
homogeneous polynomials (in $n+1$ variables) $f_1,\ldots,f_r$ of
degrees $d_1,\ldots,d_r$ respectively, the subscheme
$X=V(f_1,\ldots,f_r)\subset\P^n_k$ defined by these equations contains
a linear space $\P^m_k$ of dimension $m$.

We also define $n=l(d_1,\ldots,d_r;m;k)$ to be the smallest integer
(possibly infinite!) such that every subscheme $X\subset\P^n_k$ as
above, contains a linear space $P\cong\P^m_k$ and satisfies
$\CH_m(X)=\Z\cdot P$. Note that in this situation, the $m$-fold
intersection with the hyperplane class will induce an isomorphism
$\CH_m(X)\isom\Z$. Hence, if $P'\cong\P^m_k$ is any other linear
subspace contained in $X$ then we also have $\CH_m(X)=\Z\cdot P'$.

We now prove various inequalities about these numbers. The first and
trivial inequality is
$$ n(d_1,\ldots,d_r;m;k) \leq
        n(d_1,\ldots,d_s;n(d_{s+1},\ldots,d_r;m;k);k). \tag4.2.1$$
for any $s$ between $1$ and $r$ such that $n(d_{s+1},\ldots,d_r;m;k)$
is finite.

\proclaim{Lemma 4.2.2} We have the inequality
$$ n(d;m;k) \leq \max\{n(d;0;k), n(1,\ldots,d;m-1;k)\}. $$
\endproclaim

\demo{Proof} Let us choose $n$ greater than or equal to the number on
the right hand side of the inequality. Then we have an $k$-rational
point $p$ on $X$. We blow up the point $p$ and project to obtain
$f:\tilde{X}\to\P^{n-1}$ and $g:\tilde{\P^n}\to\P^{n-1}$.  The
condition that the fibre of $f$ at a point of $\P^{n-1}$ is the same
as the fibre of $g$ at this point gives us one equation of each degree
between 1 and $d$ on $\P^{n-1}$ (these are the {\it polars} of the
equation defining $X$, with respect to the point $p$). Let $Y$ be the
zero locus of these equations. Since $n\geq n(1,\ldots,d;m-1;k)$ we
have a linear space of dimension $m-1$ in $Y$. The inverse image in
$\tilde{X}$ pushes down to a linear subspace of dimension $m$ in $X$.
\qed\enddemo

We now observe that the inequalities (4.2.1) and (4.2.2) give us a
weaker version of a result of Leep and Schmidt \cite{12}:

\proclaim{Proposition 4.2.3}\rom{(Leep-Schmidt)} Let $k$ be a field
such that $n(d_r;0;k)$ is finite. Then $n(d_1,\ldots,d_r;m;k)$ is
finite for all integers $m$. \endproclaim

We now obtain an inequality for $l(d_1,\ldots,d_r;m;k)$.

\proclaim{Lemma 4.2.4} Let us define $l$ to be the supremum of the
numbers $l(d_1-1,\ldots,d_r-1;m';k')$, where $m'$ runs over all
numbers between 0 and $m$, and $k'$ runs over all finitely generated
field extensions of $k$ which have transcendence degree $m-m'$. If $l$
is finite, then we have the inequality $$ l(d_1,\ldots,d_r;m;k) \leq
n(d_1,\ldots,d_r;l;k). $$ \endproclaim

\demo{Proof} Let $n$ be a number greater than or equal to the right
hand side of the inequality. Then we have a linear subspace $\P^l_k$
in $X=V(f_1,\ldots,f_r)$. We blow up this subspace and project to
obtain a morphism $\tilde{X}\to\P^{n-l-1}$ and also
$\tilde{\P^n}\to\P^{n-l-1}$. Let $E$ and $\P^l_k \times \P^{n-l-1}_k$
be the exceptional divisors for the blow ups. The morphisms
$\tilde{X}\to\P^{n-l-1}_k$ and $E\to\P^{n-l-1}_k$ are families of
intersections of multidegree $(d_1-1,\ldots,d_r-1)$ in the
$\P^{l+1}$-bundle (resp. $\P^l$-bundle) $\tilde{\P^n}\to\P^{n-l-1}$
(resp. $\P^l_k \times \P^{n-l-1}_k \to \P^{n-l-1}_k$).

Let $A$ be any subvariety of $\tilde{X}$ of dimension $m$. Let $B$ be
its image in $\P^{n-l-1}_k$. By induction on the dimension of $B$, we
will show that $A$ is rationally equivalent to a $m$-cycle supported
on $E$. If $\dim B = -\infty$, then $A$ is empty and hence we are done.

Let $\dim(A/B) = m'$ and $k'=k(B)$, the function field of $B$; let
$$ E_{k'} = E \times_{\P^{n-l-1}_k} \Spec k(B)
        \subset X_{k'} = \tilde{X} \times_{\P^{n-l-1}_k} \Spec k(B) $$
Since $l$ has so been chosen that
$l\geq l(d_1-1,\ldots,d_r-1;m';k')$, we have a $P_{k'}=\P^{m'}_{k'}$
contained in $E_{k'}$ such that
$$ \Z\cdot P_{k'} \onto \CH_{m'}(E_{k'}) \onto \CH(X_{k'}).\tag A$$
Let $C\subset E$ be the closure of $P_{k'}$. Then by (A) we have
an integer $a$ such that the $m$-cycle $A - a\cdot C$ on
$X_B=\tilde{X} \times_{\P^{n-l-1}_k} B$ restricts to
a rationally trivial cycle in $X_{k'}$. If $B$ has dimension zero then
we note that $X_B=X_{k'}$ so that we are done. In any case, by
(4.1.3) we see that $A - a\cdot C$ is rationally equivalent in
$X_B$ to a cycle supported over a proper subscheme of $B$. By induction
on the dimension of $B$ we are done.

Hence we have shown that the natural morphism $\CH_m(E) \onto
\CH_m(\tilde{X})$ is a surjection. Now the image of $\CH_m(E) \to
\CH_m(X)$ factors through $\CH_m(\P^l_k)$ and so we have the result.
\qed\enddemo

Combining Proposition~(4.2.3) and Lemma~(4.2.4) we have

\proclaim{Theorem 4.2.5} Let $k$ be a field such that for all finitely
generated field extensions $k'\contain k$ of transcendence degree at
most $m$ there is a uniform bound for $n(d_r;0;k')$. Then,
$l(d_1,\ldots,d_r;m;k)$ is finite.  \endproclaim

\remark{Remark 4.2.6} Note that if $F$ is an algebraically closed
field, and $k$ is finitely generated extension of $F$ then it always
satisfies the hypothesis of (4.2.3) and (4.2.5) by results of Lang and
Nagata (see \eg \cite{7}).  \endremark

\remark{Remark 4.2.7} If we only wish to obtain (4.2.5) upto torsion
where $k$ is replaced by a universal domain $K$, the bound obtained in
(4.2.4) can be improved. The same proof will show that we need only
take $l$ to be the supremum of the integers
$l(d_1-1,\ldots,d_r-1;m';K)$, where $m'$ runs over all numbers between
0 and $m$.  \endremark

\heading 5. The Method of Bloch and Srinivas \endheading

We first generalise the key Proposition of \cite{1}.

\proclaim{Proposition 5.1} Let $X$ be a smooth projective variety over
$k$ and $K$ be a universal domain over $k$. Let $V_l\subset X$ be
subschemes for each $l\in[0,m]$ such that $\CH_l((X-V_l)_K)=0$ for
each $l$. Then for each $l$ we have cycles $\Gamma_l\in\CH^{\dim
X}(X\times X)_{\Q}$ such that support of $\Gamma_l$ is contained in
$V_l\times X$ and a cycle $\Gamma^{m+1}\in\CH^{\dim X}(X\times
X)_{\Q}$ such that support of $\Gamma^{m+1}$ is contained in $X\times
W$, where $W\subset X$ is pure of codimension $m+1$, so that if
$\Delta_X\in\CH^{\dim X}(X\times X)_{\Q}$ is the class of the diagonal
we have an equation
$$ \Delta_X = \Gamma_0 + \cdots + \Gamma_m + \Gamma^{m+1}. \tag 5.2$$
\endproclaim

\demo{Proof} The result is obvious for $m=-1$ with
$\Gamma^0=\Delta_X$. We prove the general case by induction on $m$.
Assume that the result is known for $m-1$. We then have $\Gamma^m$
which has support contained in $X\times W'$ where $W'$ is pure of
codimension $m$ in $X$. Let $P$ be any geometric generic point of $W'$
with values in $K$, \ie choose a $k$-embedding of the function field
of some component of $W'$ into $K$.  The restriction of $\Gamma^m$ to
$X\times P$ gives an element of $\CH_m(X_K)$. By assumption, this
cycle is then rationally equivalent to a cycle
$\gamma^m_P\in\CH_m((V_m)_K)$. Choose one such $P$ for each component
of $W'$ and let $\Gamma_m$ be the sum over these $P$'s of the closures
in $X\times W'$ of $\gamma^m_P$. The difference
$\xi=\Gamma^m-\Gamma_m$ restricts to zero on $X\times P$ for each
choice of $P$. By (4.1.3) we see that there is a subscheme $W$ of $W'$
of pure codimension 1 and a cycle $\Gamma^{m+1}\in\CH^{\dim X}(X\times
X)$ supported on $X\times W$ such that $\Gamma^m - \Gamma_m$ is
rationally equivalent to $\Gamma^{m+1}$.  Hence we have the result.
\qed\enddemo

\proclaim{Lemma 5.3} Let $X$ be a smooth projective variety over $k$,
$V\subset X$ be a closed subscheme and $W\subset X$ be a closed
subscheme of pure codimension $m+1$. Moreover assume that for each
$l\in[0,k]$ we have cycles $\Gamma_l\in\CH^{\dim X}(X\times X)_{\Q}$
with support contained in $V\times X$, and a cycle
$\Gamma^{m+1}\in\CH^{\dim X}(X\times X)_{\Q}$ with support in $X\times
W$, and that these cycles satisfy (5.2). Then for each $l$ we have a
morphism $\CH^l(V)_{\Q}\to\CH^l(X)_{\Q}$ given by $\sum[\Gamma_l]_*$.
\roster
\item The morphism $\CH^l(V)_{\Q}\onto\CH^l(X)_{\Q}$ is a surjection
for $l\leq m$.
\item The cokernel of $\CH^{m+1}(V)_{\Q}\to\CH^{m+1}(X)_{\Q}$ is
finite dimensional.
\item The cokernel of $\CH^{m+2}(V)_{\Q}\to\CH^{m+2}(X)_{\Q}$ is
weakly representable.
\endroster
\endproclaim
\demo{Proof} By (5.2) we only need to look at the image of
$[\Gamma^{m+1}]_*$ for each $l$. We have a factoring of this via the
Gysin morphism $\CH^{l-m-1}(W)_{\Q}\to\CH^l(X)_{\Q}$. The results then
follow from the facts: (1) For $a<0$ we have $\CH^a(W)=0$, (2)
$\CH^0(W)_{\Q}$ is a finite dimensional $\Q$-vector space, and (3)
$\CH^1(W)$ is representable.  \qed\enddemo

Now in the situation where $X$ is a smooth subvariety of $\P^n$
defined by $r$ equations of degrees $d_1\leq\cdots\leq d_r$, and
$n\geq l(d_1,\ldots,d_r;m;k)$ for an algebraically closed field $k$
(for notation see (4.2.1)). We have a linear subspace $\P^m\subset X$
such that we can apply (5.1) with $V_l=\P^l\subset\P^m$ for all
$l\in[0,m]$. Then we can apply (5.3) with $V=\P^m$ to conclude
\proclaim{Theorem 5.4} Let $X$ be a smooth subvariety of $\P^n_k$
defined by $r$ equations of degrees $d_1\leq\cdots\leq d_r$ and let
$m$ be a positive integer. If $n\geq l(d_1,\ldots,d_r;m;k)$,
then
\roster
\item For $l\in[0,m]$ we have $\CH^l(X)_{\Q}=\Q$.
\item $\CH^{m+1}(X)_{\Q}$ is finite dimensional.
\item $\CH^{m+2}(X)_{\Q}$ is representable.
\endroster
\endproclaim

Further applications of (5.1) may be found by applying the methods of
\cite{1}. We mention one consequence which is to conclude the validity
of the generalised Hodge conjecture of Grothendieck (see \cite{8}) in
some cases.
\proclaim{Proposition 5.5} In the situation of (5.4) we have
$\H^l(X)=\N^{m+1}\H^l(X)$ for all $l$ such that $l\geq (2m+2)$.
\endproclaim

\demo{Proof} We note that $[\Gamma^{m+1}]_*$ is an endomorphism of
$H^l(X)$ that factors through the Gysin homomorphism
$\H^{l-2(m+1)}(\tilde{W})\to\H^l(X)$, where $\tilde{W}$ is the
desingularization of $W$. By definition, the image of this Gysin
homomorphism is in $\N^{m+1}\H^l(X)$. Further, the endomorphism
$[\Gamma_l]^*$ of $\H^l(X)$ factors through $\H^l(V)$. In our case
$V=\P^m$, so that $\H^l(V)=\N^{l/2}\H^l(V)$. \qed\enddemo

\enddocument